\documentclass[11pt,a4paper]{amsart}
\usepackage{amscd,amssymb,amsopn,amsmath,amsthm,graphics,amsfonts,enumerate,verbatim,calc}

\headsep=1cm \numberwithin{equation}{section}
\newcommand{\To}{\rightarrow}

\newcommand{\id}{{\rm{id}}}

\newcommand{\Gid}{{\rm{Gid}}}

\newcommand{\Coker}{{\rm{Coker}}}
\newcommand{\Ker}{{\rm{Ker}}}
\newcommand{\Ass}{{\rm{Ass}}}

\newcommand{\Spec}{{\rm{Spec}}}

\newcommand{\depth}{{\rm{depth}}}

\newcommand{\Hom}{{\rm{Hom}}}
\newcommand{\Ext}{{\rm{Ext}}}

\newtheorem{theorem}{Theorem}[section]
\newtheorem{corollary}[theorem]{Corollary}
\newtheorem{lemma}[theorem]{Lemma}

\theoremstyle{definition}
\newtheorem{definition}[theorem]{Definition}
\newtheorem{definitions and notations}[theorem]{Definitions and Notations}

\newtheorem{remark}[theorem]{Remark}

\theoremstyle{plain}

\theoremstyle{definition}

\numberwithin{equation}{section}

\begin{document}

\title{Gorenstein injective dimension and a generalization of Ischebeck
Formula}
\author{Reza Sazeedeh}

\address{Department of Mathematics, Urmia University, P.O.Box: 165, Urmia, Iran-And\\
School of Mathematics, Institute for Research in Fundamental
Sciences (IPM), P. O. Box: 19395-5746, Tehran, Iran}
\email{rsazeedeh@ipm.ir}



\begin{abstract}
Let $(R,\frak m)$ be a commutative Noetherian local ring and let
$M$ and $N$ be finitely generated $R$-modules of finite injective
dimension
 and finite Gorenstein injective
dimension, respectively. In this paper we prove a generalization
of Ischebeck Formula, that is
$\depth_RM+\sup\{i|\hspace{0.1cm}\Ext_R^i(M,N)\neq 0\}=\depth R.$
\end{abstract}

\maketitle


\section{Introduction}
\hspace{0.4cm} Throughout this paper, $(R,\frak m)$ is a
commutative Noetherian local ring with the maximal ideal $\frak
m$.

 Our motivation to do this work is the first steps of the solution of a conjecture of Bass given by Levin and Vasconcelos in 1968  [LV] when $R$ admits a
 finitely generated $R$-module
 of injective dimension $\leq 1$. One of the main ingredients of their proof is to use a formula which relates
 the depth of a module to the depth of the ring via a non-vanishing of Ext when $R$ admits a finitely generated $R$-module of finite
 injective dimension. More precisely, if $M$ is a finitely generated $R$-module of finite injective dimension, then for
 any finitely generated $R$-module $L$, there is the following equality
$$\depth_RL+\sup\{i|\hspace{0.1cm}\Ext_R^i(L,M)\neq 0\}=\depth R.$$

 This formula was simultaneously proved by Ischebeck [I] using an easy induction argument on $\depth_R L$.
In this paper we give and prove a generalization of this formula
existing a finitely generated $R$-module of finite Gorenstein
injective dimension which turns out to be much more complicated.
To be more precise, we show that if $M$ and $N$ are finitely
 generated $R$-modules
 of finite injective dimension and finite Gorenstein injective
dimension respectively, then we have
$$\depth_RM+\sup\{i|\hspace{0.1cm}\Ext_R^i(M,N)\neq 0\}=\depth R.$$
As an application, we show that if $R$ is Cohen-Macaulay and $N$
is a finitely generated $R$-module of
 finite Gorenstein injective dimension, then there exists a finitely generated $R$-module $M$ of finite injective dimension such that
 $\Hom_R(M,N)\neq 0$ and $\Ext_R^i(M,N)=0$ for all $i>0$.


\section{The main results}
Throughout this section, let $(R,\frak m)$ be a local ring and we
let $\hat{R}$, the completion of $R$ with respect to the maximal
ideal $\frak m$. For an $R$-module $M$, a prime ideal $\frak p$
of $R$
 and an integer $i$, the $i$-th Bass number of $M$ with respect to $\frak p$ is the (cardinal) number $\mu_R^i(\frak p,M)$,
that is the dimension of $\Ext_{R_{\frak p}}^i(R_{\frak p}/\frak
p R_{\frak p},M_{\frak p})$ as a vector-space on the residue
field $R_{\frak p}/\frak p R_{\frak p}$.

\medskip

\begin{lemma}\label{2.2}
Let $N$ be a finitely generated $R$-module of finite Gorenstein
injective dimension and let $M$ be a finitely generated
$R$-module of finite injective dimension with $\depth_R M=\depth
R=r>0$. Then we have the following equality
$\sup\{i|\hspace{0.1cm}\Ext_R^i(M,N)\neq 0\}=0.$
\end{lemma}
\begin{proof}
We note that $\id_RM=\id_{\hat{R}}(M\otimes_R\hat{R})=r$ and by
virtue of [FFr, Theorem 3.6] we have
 $\Gid_{\hat{R}}(N\otimes_R\hat{R})<\infty$. Therefore it follows from [CS, Corollary 2.3] that
$\Gid_RN=\Gid_{\hat{R}}(N\otimes_R\hat{R})=r$, and so without
loss of generality the claim, we may assume that $R$ is complete.
Since $M$ is finitely generated and $\depth M=r$, the proof of
[FFGR, Proposition 2.2] and [PS, Theorem 4.10] imply that for each
$i\geq 0$, there is the following isomorphisms
$$\Ext^i_R(M,N)\cong\Ext^{i+r}_R(E^r,N)\cong \oplus_{\mu^r_R(\frak
m,M)}\Ext^{i+r}_R(E(R/\frak m),N).$$ Now, using [EJ2, Corollary
4.4], we have the following equalities which gives the assertion
$$r=\Gid_RN=\sup\{i|\hspace{0.1cm}\Ext^i_R(E(R/\frak m),N)\neq
0\}.$$
\end{proof}

\medskip
\begin{definition}
Let $N$ be an $R$-module. A complex $\dots\To E_1\To E_0\To N\To
0$ is called an {\it injective resolvent} of $N$ if for each
injective $R$-module $E$, the functor $\Hom_R(E,-)$ leaves it
exact. Let $M$ and $N$ be two $R$-modules and let $\dots\To
E_1\To E_0\To N\To 0$ and $0\To M\To E^0\To E^1\To\dots$ be
injective resolvent and injective resolution for the
corresponding modules $N$ and $M$. Then the 3rd quadrant double
complex ($\Hom(E^i, E_j))_{i,j}$ is such that the two associated
spectral sequences collapse. This implies that we can compute
derived functors of Hom using either the injective resolution of
$M$ or the injective resolvent of $N$ when it exists. These left
derived funetors will be denoted $\Ext_n^R(M,N)$. From the
definition of $\Ext_0^R(M,N)$ it is clear that there is a natural
transformation $\Ext_0^R(M,N)\To\Ext^0_R(M,N)=\Hom_R(M,N)$. The
image of $\Ext_0^R(M,N)$ in Hom$_R(M,N)$ consists of those linear
maps $M\To N$ which can be factored through an injective
$R$-module.
 We will let $\overline{\Ext}_0^R(M,N)$ and
$\overline{\Ext}^0_R(M,N)$ denote the kernel and cokernel of the
natural transformation $\Ext_0^R(M,N)\To \Hom_R(M,N)$.
 Following [EJ1], an $R$-module $N$ is said to be {\it mock finitely
generated} if for each finitely generated $R$-module $M$ and each
$i\geq 1$ each of the modules $\Ext_R^i(M,N),\Ext_ i^R(M,N),
\overline{\Ext}_0^R(M,N)$ and $\overline{\Ext}^0_R(M,N)$ is
finitely generated.
\end{definition}
\medskip

\begin{remark}\label{2.4}
Foxby in [F, Corollary 4.5] proved that if $M$ is a non-zero
finitely generated $R$-module with $t=\depth_RM<\dim R$, then
$\mu^{t+1}_R(\frak m,M)>\mu^{t}_R(\frak m,M)>0$. In the following
lemma which we use from this fact, $M$ is of finite injcetive
dimension
 and so $R$ is Cohen-Macaulay. In this case the result has an easy proof which we mention it. In fact, without
 loss of generality we may assume that $R$ is complete and so $R$ is a homomorphic image
of a Gorensatein local ring $(A,\frak n)$ and we have
$\depth_RM=\depth_AM$. By virtue of [F, Remark 2.7] there exists
a complex of finitely generated free
 $R$-modules $\mathcal{F}:=\dots\To F_{t+1}\To F_t\To 0$ such that the rank of each $F_i$ is $\mu^{i}_R(\frak m,M)$ and $H_i(\mathcal{F})=0$ for all
$i<\depth_RM$ or $i> \dim M$ and moreover $\dim
H_i(\mathcal{F})\leq i$ for all $i$. By using the local duality
theorem for local cohomology and this fact that $H_{\frak
n}^t(M)\neq 0$
 we can deduce that $H_t(\mathcal{F})\neq 0$. Now consider the exact sequence $0\To Z\To F_{t+1}\To F_t\To H_t(\mathcal{F})\To 0$.
As $R$ is Cohen-Macaulay, we have $\dim Z=\dim R$ and since
$t<\dim R$, there exists a $\frak p\in\Ass Z$ such that $\dim
R/\frak p>t$ and so
 $H_t(\mathcal{F})_{\frak p}=0$. Now if we localize the above exact sequence to $\frak p$ we deduce that rank$F_{t+1}>{\rm rank}F_t$.
 Keeping in mind this fact we have the following lemma.
\end{remark}

\medskip

\begin{lemma}\label{2.5}
Let $N$ be a finitely generated $R$-module of finite Gorenstein
injective dimension and let $M$ be a finitely generated
$R$-module of finite injective dimension with $\depth M=0$ and
$\depth R=1$, then
 we have the equality $\sup\{i|\hspace{0.1cm}\Ext_R^i(M,N)\neq 0\}=1.$
\end{lemma}
\begin{proof}
Similar to the argument mentioned in the previous lemma we may
assume that $R$ is complete. As $\Gid_RN=\depth R=1$, it is clear
that $\Ext_R^i(M,N)=0$ for all $i>1$ and the same reasoning
mentioned in Lemma  \ref{2.2}  and [EJ2, Corollary 4.4] imply that
$\Ext_R^1(E(R/\frak m),N)\neq 0$ and $\Ext_R^1(E(R/\frak p),N)=0$
for all $\frak p\in\Spec R\setminus\{\frak m\}$. Now, it remains
to show that $\Ext_R^1(M,N)\neq 0$. Assume that
$\Ext_R^1(M,N)=0$. As $\Gid_RN=\id_RM=1$, there exists an exact
sequence of $R$-modules $0\To N\To G\To K\To 0$ such that $G$ is
Gorenstein injective and $K$ is injective and moreover there is a
minimal injective resolution for $M$ as $0\To M\To E^0\To E^1\To
0.$ we note that since $\depth M=0$, the injective module
$E(R/\frak m)$ appears in each $E^i$. Now in view of [EJ1,
Proposition 2.4], we have the following commutative diagram with
the rows and columns exact
\begin{center}\setlength{\unitlength}{3cm}
\begin{picture}(4.6,1.9)

 \put(.9,0){\makebox(0,0){$0$}}

 \put(0.85,0.3){\makebox(0,0){$\overline{\Ext}^0_R(M,N)$}}
  \put(2,0.3){\makebox(0,0){$0$}}
 \put(3.1,0.3){\makebox(0,0){$0$}}

\put(0.895,.22){\vector(0,-1){0.15}}

\put(-0.15,0.7){\makebox(0,0){$0$}}
\put(0.858,0.7){\makebox(0,0){$\Hom_R(M,N)$}}
\put(1.93,0.7){\makebox(0,0){$\Hom_R(M,G)$}}
\put(3.03,0.7){\makebox(0,0){$\Hom_R(M,K))$}}
\put(4.19,0.7){\makebox(0,0){$\Ext_R^1(M,N)=0$}}

 \put(-0.05,.7){\vector(1,0){0.4}}
\put(1.25,.7){\vector(1,0){0.3}}
 \put(2.33,.7){\vector(1,0){0.25}}
\put(3.5,.7){\vector(1,0){0.172}}

\put(0.895,.6){\vector(0,-1){0.213}}
 \put(1.995,.6){\vector(0,-1){.213}}
\put(3.1,.6){\vector(0,-1){0.213}}

\put(-0.15,1.1){\makebox(0,0){$0=\Ext_1^R(M,K)$}}
\put(0.85,1.1){\makebox(0,0){$\Ext_0^R(M,N)$}}
\put(1.9,1.1){\makebox(0,0){$\Ext_0^R(M,G)$}}
\put(3,1.1){\makebox(0,0){$\Ext_0^R(M,K)$}}
\put(4,1.1){\makebox(0,0){$0$}}

\put(0.33,1.1){\vector(1,0){0.15}}
 \put(1.241,1.1){\vector(1,0){0.3}}
 \put(2.27,1.1){\vector(1,0){0.3}}
\put(3.45,1.1){\vector(1,0){0.475}}

\put(0.895,1){\vector(0,-1){0.213}}
 \put(1.995,1){\vector(0,-1){.213}}
\put(3.1,1){\vector(0,-1){0.213}}
 \put(3.15,.89){$\delta$}

\put(0.85,1.5){\makebox(0,0){$\overline{\Ext}_0^R(M,N)$}}
\put(1.9,1.5){\makebox(0,0){$\overline{\Ext}_0^R(M,G)$}}
\put(3,1.5){\makebox(0,0){$\overline{\Ext}_0^R(M,K)$}}

\put(1.241,1.5){\vector(1,0){0.283}}
 \put(2.253,1.5){\vector(1,0){0.293}}

\put(0.895,1.4){\vector(0,-1){0.213}}
 \put(1.995,1.4){\vector(0,-1){.213}}
\put(3.1,1.4){\vector(0,-1){0.213}}
 \put(0.92,.9){$\beta$}
\put(2.03,.9){$\alpha$}

\put(0.9,1.8){\makebox(0,0){$0$}}
\put(1.99,1.8){\makebox(0,0){$0$}}
 \put(3.083,1.8){\makebox(0,0){$0$}}

\put(0.895,1.735){\vector(0,-1){0.1679}}
 \put(1.958,1.625){$\|$}
\put(3.037,1.625){$\|$}
\end{picture}
\end{center}
The Snack Lemma implies the following exact sequence of
$R$-modules $$0=\Ker\delta\To\Coker\beta\cong
\overline{\Ext}_R^0(M,N)\To \Coker\alpha=0$$ which implies that
$\overline{\Ext}_R^0(M,N)=0$. On the other hand, by virtue of
[EJ1, Proposition 1.6] there is the following exact sequence
$$0=\overline{\Ext}_R^0(M,N)\To
\Ext_R^1(E^1,N)\To\Ext_R^1(E^0,N)\To\Ext_R^1(M,N)=0$$ and so
since $M$ is finitely generated, we have the following
isomorphisms
$$\oplus_{\mu^1_R(\frak m,M)}\Ext_R^1(E(R/\frak
m),N)\cong\Ext_R^1(E^1,N)\cong\Ext_R^1(E^0,N)\cong\oplus_{\mu^0_R(\frak
m,M)}\Ext_R^1(E(R/\frak m),N).$$ The preceding paragraph implies
that $\mu^1_R(\frak m,M)=\mu^0_R(\frak m,M)$. But as $\depth_R
M<\depth R\leq \dim R$, the last statement is in contradiction
with Remark \ref{2.4} and so $\Ext_R^1(M,N)\neq 0$.
\end{proof}

\medskip
\begin{theorem}
Let $M$ be a finitely generated $R$-module of finite injective
dimension and let $N$ be a finitely generated $R$-module of
finite Gorenstein injective dimension. Then there is the
following equality
 $$\depth_RM+\sup\{i|\hspace{0.1cm}\Ext_R^i(M,N)\neq 0\}=\depth R.$$
\end{theorem}
\begin{proof}
As all invariants in the assumption and the assertion are
well-behaviour with respect to the completion of $R$ in $\frak
m$, we may assume that $R$ is complete. At first assume that
$\depth R=0$. In this case $\id_RM=0=\Gid_RN$, and so
$M=\oplus_{\mu^0_R(\frak m,M)}E(R/\frak m)$. Therefore the result
follows by the same reasoning in Lemma \ref{2.2}  and using [EJ2,
Corollary 4.4] in this case. Assume that $\depth R=r>0$ and
$\depth _RM=s$ and we also note that $s\leq \dim M\leq \id_R
M=r$. If $r=s$, then the assertion follows by Lemma \ref{2.2} .
Therefore we assume that $s<r$. In this case if $r=1$, then the
assertion follows by Lemma \ref{2.5}. Therefore we may assume that
$r\geq 2$. As $\Gid_RN=r$, there is an exact sequence of
$R$-modules $$0\To N\To G\To K\To 0\hspace{0.2cm}(\dag)$$ such
that $G$ is mock finitely generated Gorenstein injective and $K$
is mock
 finitely generated with $\id_RK=r-1$ (see the proof of [S, Proposition 5.3]). We claim that
 $\mu^{r-1}_R(\frak m,K)>0$, otherwise there exists a prime ideal
$\frak p$ contained in $\frak m$ such that $\Ext_R^{r-1}(R/\frak
p,K)\neq 0$.
 Now, since $\dim R/\frak p=t>0$, a similar proof mentioned in [B,
Lemma 3.1] implies that $\mu^{r-1+t}_R(\frak m,K)\neq 0$, that is
in contradiction with id$_RK=r-1$.
 Now, we proceed the rest of proof by induction on $s$. If $s=0$, there is a
monomorphism $0\To R/\frak m\To M$ and so there is an epimorphism
$\Ext^{r-1}_R(M,K)\To \Ext^{r-1}_R(R/\frak m,K)\To 0$. Since
$\Ext^{r-1}_R(R/\frak m,K)\neq 0$, we have $\Ext^{r-1}_R(M,K)\neq
0$. On the other hand it is clear that $\Ext^i_R(M,K)=0$ for all
$i>r-1$ as id$_RK=r-1$. Now, if we apply the functor
$\Hom_R(M,-)$ to that exact sequence $(\dag)$, the assertion
follows easily in this case. Suppose that $s>0$ and the assertion
has been proved for all values smaller than $s$ and so we prove
it for $s$. In this case there exists an element $x\in\frak m$
which is an $M$-regular. Hence there exists the following exact
sequence of $R$-modules $$0\To M\stackrel{x.}\To M\To M/xM\To
0\hspace{0.2cm}(\ddag).$$ Applying the functor $\Hom_R(-,N)$ to
the exact sequence $(\ddag)$, for each $i$, induces the following
exact sequence $$\Ext_R^i(M,N)\stackrel{x.}\To \Ext_R^i(M,N)\To
\Ext_R^{i+1}(M/xM,N).$$ For each $i>r-s$, we have $i+1>r-s+1$,
and so the induction hypothesis on $M/xM$ implies that
$\Ext_R^{i+1}(M/xM,N)=0$ as $\depth_RM/xM=s-1$. Hence Nakayama's
Lemma implies that $\Ext_R^i(M,N)=0$ for all $i>r-s$. On the
other hand the following exact sequence
$$\Ext_R^{r-s}(M,N)\stackrel{x.}\To \Ext_R^{r-s}(M,N)\To
\Ext_R^{r-s+1}(M/xM,N)\To 0$$  and the induction hypothesis on
$M/xM$ imply that $\Ext_R^{r-s}(M,N)\neq 0$.
\end{proof}

\medskip
\begin{corollary}
Let $(R,\frak m)$ be a Cohen-Macaulay local ring and let $N$ be a
finitely generated $R$-module of finite Gorenstein injective
dimension.
 Then there exists a finitely generated $R$-module $M$ of finite injective dimension such that $\Hom_R(M,N)\neq 0$ and
$\Ext_R^i(M,N)=0$ for all $i>0$.
\end{corollary}
\begin{proof}
Let $\depth R=\dim R=d$. Then there exist a system of parameter
 $x_1,\dots, x_d\in\frak m\setminus \frak m^ 2$. Consider $R_i=R/(x_1,\dots,x_i)R$ and we know that
  $R/(x_1,\dots,x_d)$ is Artinian and injective envelope of $R/\frak m$, say $M$, is a finitely generated $R_d$-module.
By virtue of [LV, Theorem 3.1], we have id$_{R_{d-1}}M=1$ and
repeating this way we have id$_RM=d$ and $M$ is a finitely
generated $R$-module. We note that by constructing this module
and the same reasoning mentioned in [LV, Remarks. p. 319]
 we have $\mu^{d-i}_{R_i}(\frak m R_i,M)=1$. We claim that $\depth_RM=d$. Otherwise, assume that $\depth _R M=t<t+1=\depth {R_{d-t-1}}\leq d$.
 As $R_{d-1}$ is a finitely generated $R$-module, by virtue of [LV, Lamma, p. 317], we have $\depth_{R_{d-t-1}}M=\depth_RM=t$.
 On the other hand, Remark \ref{2.4} implies that $1=\mu_{R_{d-t-1}}^{t+1}(\frak m R_{d-t-1},M)>\mu_{R_{d-t-1}}^{t}(\frak m R_{d-t-1},M)>0$ which is
impossible and so $\depth_RM=d$. Now, Theorem 2.5 implies that
$\Hom_R(M,N)\neq 0$ and $\Ext_R^ i(M,N)=0$ for all $i>0$.
\end{proof}

\end{document}